\documentclass[12pt]{article}
\title{Symplectic Resolutions for Coverings of Nilpotent Orbits}
\author{Baohua FU}
\usepackage{amsmath,amssymb,amsthm,amscd}
\chardef\bslash=`\\
\newtheorem{Thm}{Theorem}
\newtheorem{Cor}{Corollary}
\newtheorem{Lem}{Lemma}
\newtheorem{Prop}{Proposition}
\newtheorem{Def}{Definition}
\newtheorem{Rque}{Remark}

\newtheorem{Conj}{Conjecture}

\def\cit{{\mathbb C}}

\def\zit{{\mathbb Z}}
\def\0{{\mathcal O}}
\def\g{{\mathfrak g}}

\begin{document}
\maketitle

\begin{abstract}
Let $\0$ be a nilpotent orbit in a semisimple complex Lie algebra $\g$. Denote by $G$ the simply connected Lie group with Lie algebra $\g$.
For a $G$-homogeneous covering $M \rightarrow \0$, let $X$ be the normalization of $\overline{\0}$ in the function field of $M$. In this note,
we study the existence of symplectic resolutions for such coverings $X$.
\end{abstract}

\section{Introduction}
Let $G$ be a simply connected semi-simple complex Lie group with Lie algebra $\g$. Let $\0$ be a nilpotent adjoint orbit in $\g$. 
The fundamental group $\pi_1(\0)$ of $\0$ is finite and in general non-zero (see corollary 6.1.6 and section 8.4 of \cite{CM}). Let $p: M \rightarrow 
\0$ be a $G$-homogeneous covering of degree $d$. We denote by  $X$ the normalization of $\overline{\0}$ in the function field of $M$, 
which will be called a {\em covering} of the nilpotent orbit $\overline{\0}$.
Then proposition 1.2 of \cite{BK} says that $X$ contains $M$ as a Zariski open subset and the map $p$ extends to a finite surjective 
$G$-equivariant morphism $\bar{p}: X \rightarrow \overline{\0}$. Furthermore $G$ has finitely many orbits on $X$ and $X - M$ is of
codimension at least 2 in $X$.  Recall that a regular 2-form on a smooth algebraic variety is {\em symplectic} if it is closed
 and non-degenerate at every point.
 The Kostant-Kirillov form on the orbit $\0$ gives a symplectic form $\omega$ on $M$. 
\begin{Def}
A {\em symplectic resolution} for $X$ is a projective resolution of singularities $\pi: Z \rightarrow X$ such that the 2-form $\pi^*(\omega)$
defined a priori on $\pi^{-1}(M)$ extends to a symplectic form $\Omega$ on the whole of $Z$.
\end{Def}

As shown in \cite{Fu}, a resolution is symplectic if and only if it is
crepant.  The purpose of this note is to consider the existence of
symplectic resolutions, i.e. for which coverings of a nilpotent
adjoint orbit there exists a symplectic resolution.  It turns out that
this is a difficult problem. The following proposition gives some
examples of symplectic resolutions.
\begin{Prop}
Let $P$ be a parabolic subgroup of $G$ and $M$ the unique open
$G$-orbit in $T^*(G/P)$. Let $\0$ be the orbit of a Richardson
element in the nilradical $\mathfrak{u}$ of the Lie algebra of $P$.
Then the $G$-equivariant desingularization $T^*(G/P) \rightarrow X$
is a symplectic resolution for $X$.
\end{Prop}

For the proof, see proposition 7.4 of \cite{BK}. It should be pointed
out that there exist some symplectic resolutions for coverings of
nilpotent orbits which are not of the above form. An example is the
double covering of the minimal nilpotent orbit in
$\mathfrak{sp}(2n,\cit)$ which is $X = \cit^{2n} \rightarrow
\overline{\0}_{min}$. However this is the only example where $X$ is
non-singular (see theorem 4.6 \cite{BK}).  The purpose of this note is
to prove the converse of this proposition under some additional hypotheses.

Recall that $\0$ and $\overline{\0}$ are both stable under the scaling 
action of $\cit^*$ on $\g$. This induces the Euler action of $\cit^*$ on 
$\0$ and $\overline{\0}$. Unfortunately this action does not lift to $X$ 
in general. An example is the double covering of the minimal
nilpotent orbit $\0_{min}$ in $\mathfrak{sp}(2n,\cit)$. However 
we have (see lemma 1.3 and proposition 1.4  \cite{BK})
\begin{Prop}
For any covering $X$ of $\overline{\0}$, there exists  a $\cit^*$-action 
on $X$ which lifts the square of the Euler
 action of $\cit^*$ on $\overline{\0}$,
i.e. for any $\lambda \in \cit^*$ and $x \in X$, $\bar{p}(\lambda \cdot x) = \lambda^2 \bar{p}(x)$. Furthermore there exists a unique point 
$o \in X$ such that $\bar{p}(o) = 0$. This $o$ is the unique $G$-fixed point in $X$ and also the unique $\cit^*$-fixed point in $X$.
\end{Prop}

For this $\cit^*$-action, we have $\lambda^* \omega = \lambda^2 \omega$ 
for any $\lambda \in \cit^*$, which is different to our situation
in \cite{Fu}, where the $\cit^*$-action satisfies 
$\lambda^* \omega = \lambda \omega$. This makes the situation more complicated
here. However, under the hypothesis that the degree of the covering is odd, 
we can prove the following 
\begin{Thm}\label{thm}
Let $\g$ be a semisimple complex Lie algebra and $X$ a covering of a
nilpotent adjoint orbit $\0$ in $\g$. Suppose that the degree $d$ of
the covering $\bar{p}: X \rightarrow \overline{\0}$ is odd. Then for
any symplectic resolution $\pi: Z \rightarrow X$, $Z$ is isomorphic to
$T^*(G/P)$ for some parabolic subgroup $P$ of $G$. Furthermore, under
this isomorphism, the map $\bar{p} \circ \pi$ becomes
$$ T^*(G/P) \simeq G \times^P \mathfrak{u} \rightarrow \g, \qquad 
(g,u) \mapsto Ad(g)u, $$  
where $\mathfrak{u}$ is the nilradical of $\mathfrak{p} = Lie(P)$.
\end{Thm}

This theorem generalizes our main theorem in \cite{Fu}, where we 
considered the case $d = 1$.

\section{Outline of the proof}
\begin{Lem}
The Euler action of $\cit^*$ on $\overline{\0}$ lifts to $X$ and for 
this action $\lambda^* \omega = \lambda \omega$. 
\end{Lem}
\begin{proof}
Let $R$ be the regular functions ring of $X$.  For $k \in \zit$, let
$R[k] = \{ \phi \in R | \lambda \cdot \phi = \lambda^k \phi, \lambda
\in \cit^* \}$.  Proposition 1.4 of \cite{BK} implies that if the
degree of the covering $X \rightarrow \overline{\0}$ is odd, then
$R[k] = 0$ for $k$ odd, i.e. $R = \sum_0^\infty R[2k]$.  This gives
that the Euler action of $\cit^*$ (not only its square!) on
$\overline{\0}$ lifts to a $\cit^*$-action on $X$.  For this action,
we have $\lambda^* \omega = \lambda \omega$.
\end{proof}
\begin{Lem}
The action of $G \times \cit^*$ on $X$ lifts to $Z$ such that 
$\pi: Z \rightarrow X$ and $\bar{p}\circ \pi: Z \rightarrow \overline{\0}$
is $G \times \cit^*$-equivariant. For the $\cit^*$-action 
on $Z$, we have  $\lambda^* \Omega = \lambda \Omega$
for any $\lambda \in \cit^*$.
\end{Lem}

Note that the action of $G$ (resp. $\cit^*$) on $M$ lifts to
$\pi^{-1}(M) \subset X$.  Now the proof goes along the same line as
proposition 3.1 of \cite{Fu}. Using the two lemmas, we can apply our
analysis in \cite{Fu} to complete the proof of theorem 3. For reader's
convenience, we give an outline of our method.

Let $Z^{\cit^*}$ be the fixed points subvariety in $Z$. Since
$\bar{p}\circ \pi$ is proper and the $\cit^*$-action on
$\overline{\0}$ extends to a $\cit$-action, the valuative criterion of
properness shows that there exists an attraction $q: Z \rightarrow
Z^{\cit^*}$.  As show by lemma 3.5 \cite{Fu}, the $\cit^*$-action on
the smooth connected component $Z_0 \subset Z^{\cit^*}$ containing
$q(\pi^{-1}(M))$ is definite, i.e. for $z \in Z_0, T_z^kZ = \{v \in
T_zZ | \lambda_*v = \lambda^k v \}$ is zero if $k < 0$.  Now the
equation $\lambda^* \Omega = \lambda \Omega$ implies a duality between
$T_z^kZ$ and $T_z^{1-k}Z$. Thus for $z \in Z_0$, we have $T_z Z_0
\simeq T_z^0Z \simeq T_z^1Z$ and $T_z^k Z = 0$ if $k \neq 0,1$.  This
shows that $Z_0$ is Lagrangian in $Z$. By a classical result
of Bialynicki-Birula, the attraction
$q:q^{-1}(Z_0) \rightarrow Z_0$ is $\cit^*$-equivariantly isomorphic
to $T^*Z_0 \rightarrow Z_0$.

Now we study the $G$-action on $Z_0$, which has an open-dense orbit,
namely $q(\pi^{-1}(M))$. In fact this orbit is the whole of $Z_0$.
Note that $Z_0 \subset Z^{\cit^*} \subset \pi^{-1}(o)$ is projective,
so $Z_0 = G/P$ for some parabolic subgroup $P$ of $G$.  To complete
the proof, we need to show that $Z_0 = Z^{\cit^*}$. Note that
$\pi^{-1}(o)$ is connected, so we need only to show that $Z_0$ is a
connected component of $\pi^{-1}(o)$ or equivalently a connected
component of $(\bar{p} \circ \pi)^{-1}(0)$, this is proved by using the
explicit formula for the map $T^*(G/P) \rightarrow \overline{\0}, \,
(g,u) \mapsto Ad(g)u$.

\section{Some corollaries and examples}

\begin{Cor}
Let $X$ be a covering of odd degree of a nilpotent orbit $\overline{\0}$. 
If $\0$ is not a Richardson orbit, then $X$ admits no symplectic resolution.
\end{Cor}

This follows directly from the above theorem and proposition 2.
Consider the nilpotent orbits $2A_2, 2A_2+A_1, A_5, E_6(a_1), E_6$ in
the exceptional Lie algebra of type $E_6$, here we use notations from
section 8.4 of \cite{CM}. All these nilpotent orbits have $\zit/3$ as
fundamental group, but none of them is a Richardson orbit. The above
corollary implies that the three-fold covering of any of the above
orbits admits no symplectic resolution.

\begin{Cor}
A covering $X$ of odd degree $d > 1$  of a nilpotent orbit $\0$ in 
$\mathfrak{sl}(n,\cit)$ admits no symplectic resolution.
\end{Cor}

This follows from the above theorem and theorem 3.3 \cite{Hes}, which
says that for any polarization $P$ of $\0 \subset
\mathfrak{sl}(n,\cit)$, the morphism $T^*(G/P) \rightarrow
\overline{\0}$ is birational. Recall that every nilpotent orbit in
$\mathfrak{sl}(n,\cit)$ admits a symplectic resolution (\cite{Fu}),
however the above corollary shows that their odd higher degree
coverings admit no symplectic resolution.

As an example, let $\0$ be the principal nilpotent orbit in
$\mathfrak{sl}(3,\cit)$, then $\pi_1(\0) = \zit/3$.  Let $X$ be the
three-fold covering of $\overline{\0}$. Then by the above corollary,
$X$ does not admit any symplectic resolution.  In fact, $X$ is a
so-called shared orbit in \cite{BK}, more precisely $X$ is nothing but
the closure of the minimal nilpotent $\overline{\0}_{min}$ in the
exceptional Lie algebra $\g_2$, which is known to admit no symplectic
resolution \cite{Fu}.
  
In the list of \cite{BK}, shared orbits do not admit any symplectic
resolution except in the following case: let $\0$ be the nilpotent orbit
in $\mathfrak{sp}(2n,\cit)$ corresponding to the partition
$[2,2,1,\cdots,1]$. Then $\pi_1(\0) = \zit/2$.  By results of
\cite{Fu}, the nilpotent orbit $\overline{\0}$ does not admit any
symplectic resolution.  Let $X$ be the double covering of
$\overline{\0}$. Then $X$ is isomorphic to $\overline{\0}_{min}$ in
$\mathfrak{sl}(2n,\cit)$, which admits a
symplectic resolution.

As shown by the above examples, there exist some nilpotent orbits
which admit some symplectic resolutions, but not their coverings, and
there exist some nilpotent orbits which do not admit any symplectic
resolution, while some of their coverings do admit some symplectic
resolutions.  This indicates that the problem to determine which
covering admits a symplectic resolution might be difficult and
interesting.

\section{A conjecture}
From now on, we suppose that $\g$ is a simple complex Lie algebra. Let
$\0$ be a nilpotent adjoint orbit in $\g$ and $X$ a covering of
$\overline{\0}$. Recall that $R[k] = \{ \phi \in R | \lambda \cdot
\phi = \lambda^k \phi, \lambda \in \cit^* \}$.  Let $\g' = R[2]$,
which is also a simple Lie algebra (see theorem 4.2 of \cite{BK}). Let
$G'$ be the simply connected Lie group with Lie algebra $\g'$, which
is the maximal connected Lie group of holomorphic symplectic
automorphisms of $X$ (see corollary 7.3 \cite{BK}).  By theorem 3.1
\cite{BK}, $X$ contains a Zariski open subset $M'$ such that $M
\subset M'$, and there exists a finite $G'$-covering $M' \rightarrow
\0'$ to a nilpotent orbit $\0'$ in $\g'$. We propose the following
conjecture on symplectic resolutions for $X$:
\begin{Conj}
Let $\bar{p}: X \rightarrow \overline{\0}$ be a covering of a nilpotent 
orbit in a simple Lie algebra $\g$ such that $R[1] = 0$.
 Suppose that $\pi: Z \rightarrow X$ 
is a symplectic resolution. Then there exists some parabolic 
sub-group $P'$ of $G'$ such that $Z$ is isomorphic to $T^*(G'/P')$ 
and the map $\pi$ becomes  $$T^*(G'/P') \rightarrow X. $$
\end{Conj}

\begin{Rque}
 It is proved in \cite{BK} (theorem 4.6) that $R[1] \neq 0$ if and only 
if $X$ is the double covering of $\overline{\0}_{min} \subset 
\mathfrak{sp}(2n,\cit)$. In this case, $X = \cit^{2n}$.
\end{Rque}
\begin{Rque}
By theorem 7.5 \cite{BK}, if $T^*(G/P) \rightarrow X$ is a symplectic 
resolution for $X$, then there exists some parabolic subgroup $P'$ of $G'$
such that $G/P = G'/P'$.
\end{Rque}
\begin{Rque}
If the covering $\bar{p}: X \rightarrow \overline{\0}$ is of odd degree, then by 
our theorem \ref{thm}, the conjecture is true. 
\end{Rque}

{\bf Acknowledgements:} I want to thank A. Beauville and M. Brion for helpful
discussions and suggestions. I am grateful to the referee for the  careful reading and pertinent remarks.

\end{document}